\documentclass{article}
\usepackage{amsmath}
\usepackage{amsfonts}
\usepackage{graphicx}

\newtheorem{theorem}{Theorem}

\newtheorem{proposition}[theorem]{Proposition}

\begin{document}

\title{Stochastic solutions and singular partial differential equations }
\author{R. Vilela Mendes\thanks{%
rvilela.mendes@gmail.com, rvmendes@fc.ul.pt,
https://label2.tecnico.ulisboa.pt/vilela/} \\
CMAFcIO, Faculdade de Ci\^{e}ncias, Universidade de Lisboa}
\date{ }
\maketitle

\begin{abstract}
The technique of stochastic solutions, previously used for deterministic
equations, is here proposed as a solution method for partial differential
equations driven by distribution-valued noises.
\end{abstract}

\section{Introduction}

Singular partial differential equations (SPDE's) are partial differential
equations (PDE's) which are singular in the sense that, for example,
nonlinear functions of the driving terms or of the solution itself are
ill-defined. As a result, even the question of what means to be a solution
is a non-trivial matter. Examples of current interest are:%
\begin{equation}
\partial _{t}h=\partial _{x}^{2}h+\lambda \left( \partial _{x}h\right)
^{2}+\xi ,  \label{1.1}
\end{equation}%
the Kardar-Parisi-Zhang (KPZ) model \cite{KPZ} for interface propagation in $%
d=1$;%
\begin{equation}
\partial _{t}\Phi =\Delta \Phi -\Phi ^{3}+\xi ,  \label{1.2}
\end{equation}%
the stochastic quantization equation for Euclidean $\Phi ^{4}$ field theory 
\cite{Parisi} or a model for near mean field phase transition dynamics in
dimension $d=2$ or $3$ \cite{Presutti}, $\xi $ being a space-time white
noise and a few others \cite{Carmona} \cite{Hairer3}.

Past methods to deal with these equations, and in general with nonlinear
partial differential equations driven by singular stochastic processes, are
based either on controlled rough paths \cite{Lyons} \cite{Gubinelli1},
paraproducts \cite{Gubinelli2} or a theory of regularity structures \cite%
{Hairer3} \cite{Hairer1} \cite{Hairer2}.

One reason why handling these SPDE's is a delicate matter is because, in an
iterative process of solution, nonlinear terms, products and powers of
distributions do appear and there is no canonical way to multiply arbitrary
distributions. For example in the KPZ equation, if the solution has the same
regularity structure as Brownian motion, the nonlinear term would involve
powers of the white noise, that is, products of distributions at the same
spacetime points. However, the distributions that appear in the SPDE's
solutions are not arbitrary, being constrained by the structure of the
equation itself, in particular by its linear part. Controlling this effect
is at the basis of the previously proposed solution methods. Here, in the 
\textit{stochastic solution method} that is proposed, no attempt is made to
address the general problem of multiplying arbitrary distributions. Instead
one uses the fact that, as already pointed out in \cite{Vilela-ultra}, when
one uses a stochastic process to represent the solution of a nonlinear
differential equation, the products in the nonlinear terms are traded off by
branchings of the process. Multiplication of distribution-valued quantities
may still appear when computing the functionals or the boundary measures
that represent the solution. However, at least when the distributions have
point support or when they correspond to processes that are uncorrelated at
distinct points, it turns out that, with probability one, the computation of
the process solution becomes well defined.

But, what is a stochastic solution?

The solutions of linear elliptic and parabolic equations, both with Cauchy
and Dirichlet boundary conditions, have a probabilistic interpretation.
These are classical results which may be traced back to the work of Courant,
Friedrichs and Lewy \cite{Courant} in the 1920's and became a standard tool
in potential theory \cite{Getoor} \cite{Bass1} \cite{Bass2}. For example,
for the heat equation

\begin{equation}
\partial _{t}u(t,x)=\frac{1}{2}\frac{\partial ^{2}}{\partial x^{2}}%
u(t,x)\qquad \mathnormal{with}\qquad u(0,x)=f(x)  \label{1.7}
\end{equation}%
the solution may be written either as 
\begin{equation}
u\left( t,x\right) =\frac{1}{2\sqrt{\pi }}\int \frac{1}{\sqrt{t}}\exp \left(
-\frac{\left( x-y\right) ^{2}}{4t}\right) f\left( y\right) dy  \label{1.8}
\end{equation}%
or as 
\begin{equation}
u(t,x)=\mathbb{E}_{x}f(X_{t})  \label{1.9}
\end{equation}%
$\mathbb{E}_{x}$ meaning the expectation value, starting from $x$, of the
process%
\begin{equation}
dX_{t}=dB_{t}  \label{1.10}
\end{equation}%
$B_{t}$ being the Wiener process.

Eq.(\ref{1.7}) is a \textit{specification} of the problem whereas (\ref{1.8}%
) and (\ref{1.9}) are \textit{exact} \textit{solutions} in the sense that,
given an initial condition, they both provide algorithmic means to construct
a function satisfying the specification. An important condition for (\ref%
{1.8}) and (\ref{1.9}) to be considered as solutions is the fact that the
algorithmic tools are independent of the particular solution, in the first
case an integration procedure and in the second a functional of a
solution-independent process. This should be contrasted with stochastic
processes constructed from a given particular solution, as has been done for
example for the Boltzman equation \cite{Graham}.

In contrast with the linear problems, for nonlinear partial differential
equations explicit solutions in terms of elementary functions or integrals
are only known in very particular cases. However, if a solution-independent
stochastic process is found that (for arbitrary initial conditions)
generates the solution in the sense of Eq.(\ref{1.9}), an exact stochastic
solution is obtained. In this way the set of equations for which exact
solutions are known would be considerably extended.

The stochastic solutions recently constructed for the deterministic
Navier-Stokes \cite{Jan} \cite{Waymire} \cite{Waymire2} \cite{Bhatta1} \cite%
{Orum} \cite{Ossiander}, Vlasov-Poisson, Magnetohydrodynamics and other
equations \cite{Vilela1} \cite{Vilela2} \cite{VilelaJMP} \cite{FlorianiJCP} 
\cite{Ouerdiane} define solution-independent processes for which the mean
values of some functionals are solutions to these equations. Therefore, they
are exact \textbf{stochastic solutions}.

In the stochastic solutions one deals with a process that starts from the
point where the solution is to be found, a functional being then computed
along the sample path until it reaches a boundary. For deterministic
equations the solution is an average value over all sample paths of the
process, whereas for stochastic PDE's the process itself (and its
functional) might be considered to be the solution. The localized and
parallelizable nature of the solution construction allows for the
development of efficient parallel computing algorithms \cite{Acebron1} \cite%
{Acebron2}. Provided some differentiability conditions are satisfied, the
process also handles equally well simple or very complex boundary conditions.

Stochastic solutions also provide an intuitive characterization of the
physical phenomena, relating nonlinear interactions with cascading
processes. By the study of exit times from a domain they also sometimes
provide access to quantities that cannot be obtained by perturbative methods 
\cite{VilelaZeit}.

One way to construct stochastic solutions is based on a probabilistic
interpretation of the Picard series. The differential equation is written as
an integral equation which is rearranged in such a way that the coefficients
of the successive terms in the Picard iteration obey a normalization
condition. The Picard iteration is then interpreted as an evolution and
branching process, the stochastic solution being equivalent to importance
sampling of the normalized Picard series. This method traces back its origin
to McKean's paper \cite{McKean} on the Kolmogorov-Petrovskii-Piskunov (KPP)\
equation \cite{KPP}. Hence I will refer to it as \textit{the McKean method}.

Another method is \textit{the method of superprocesses} developed by Dynkin 
\cite{Dynkin1} \cite{Dynkin2}. In this method for each equation a
superprocess on measures is generated, the solution being obtained by a
scaling process and the measure generated on the boundary. As originally
proposed the method applies only to a limited class of equations. This
limitation has however been lifted by extending from superprocesses on
measures to superprocesses on ultradistributions \cite{VilelaSPOU}.
Nevertheless, because of the limiting scaling process needed to approach the
solution, the superprocess method seems to be less efficient for the
implementation of practical solutions of PDE's. For a comparison of the two
approaches refer to \cite{Vilela-versus}.

In the past stochastic solutions have been used to find new exact solutions
of deterministic PDE's. Here one proposes this method to solve SPDE's by
treating the two examples of SPDE's listed above. The solution
representations that are obtained may in some cases coincide with those
obtained by other methods. Nevertheless they provide a very intuitive
representation of the solutions which, in addition, also provide a very
convenient framework for concrete numerical implementations.

\section{The KPZ equation}

A useful transformation for the KPZ equation (\ref{1.1}) is%
\begin{eqnarray}
Z &=&e^{-t}e^{\lambda h}  \notag \\
h &=&\frac{1}{\lambda }\left( t+\log Z\right)   \label{K1}
\end{eqnarray}%
that is, a Cole-Hopf transformation times the factor $e^{-t}$, leading to%
\begin{equation}
\partial _{t}Z=\partial _{x}^{2}Z-Z+\lambda Z\xi   \label{K2}
\end{equation}%
Of course, this transformation does not avoid the difficulties associated to
the product of distributions, because also the last term involves two
distribution-valued entities. To construct a stochastic solution the first
step is to rewrite (\ref{K2}) as an integral equation%
\begin{eqnarray}
Z\left( t,x\right)  &=&e^{-t}e^{t\partial _{x}^{2}}Z\left( 0,x\right)
+\int_{0}^{t}e^{-\left( t-s\right) }e^{\left( t-s\right) \partial
_{x}^{2}}Z\left( s,x\right) \lambda \xi \left( s,x\right)   \notag \\
&=&e^{-t}e^{t\partial _{x}^{2}}Z\left( 0,x\right)
+\int_{0}^{t}e^{-s}e^{s\partial _{x}^{2}}Z\left( t-s,x\right) \lambda \xi
\left( t-s,x\right)   \label{K3}
\end{eqnarray}%
$Z\left( t,x\right) $ being distribution-valued, the integrals in Eq.(\ref%
{K3}) involve a product of distributions and, when interpreted as usual
equalities, would be meaningless. However they become well-defined when
interpreted as a symbolic definition of an iterated stochastic process. The
product $Z\left( s,x\right) \xi \left( s,x\right) $ inside the integral
simply means that the $Z-$process samples the driving term at a branching
point. This case is the simplest instance of a branching. A term $%
Z^{3}\left( s,x\right) $, for example, would mean a branching at $\left(
s,x\right) $ of the process into three similar processes. It is in this
sense that, as stated before, undefined multiplications of distributions are
traded off by well defined branchings. 

\begin{figure}[htb]
\centering
\includegraphics[width=0.8\textwidth]{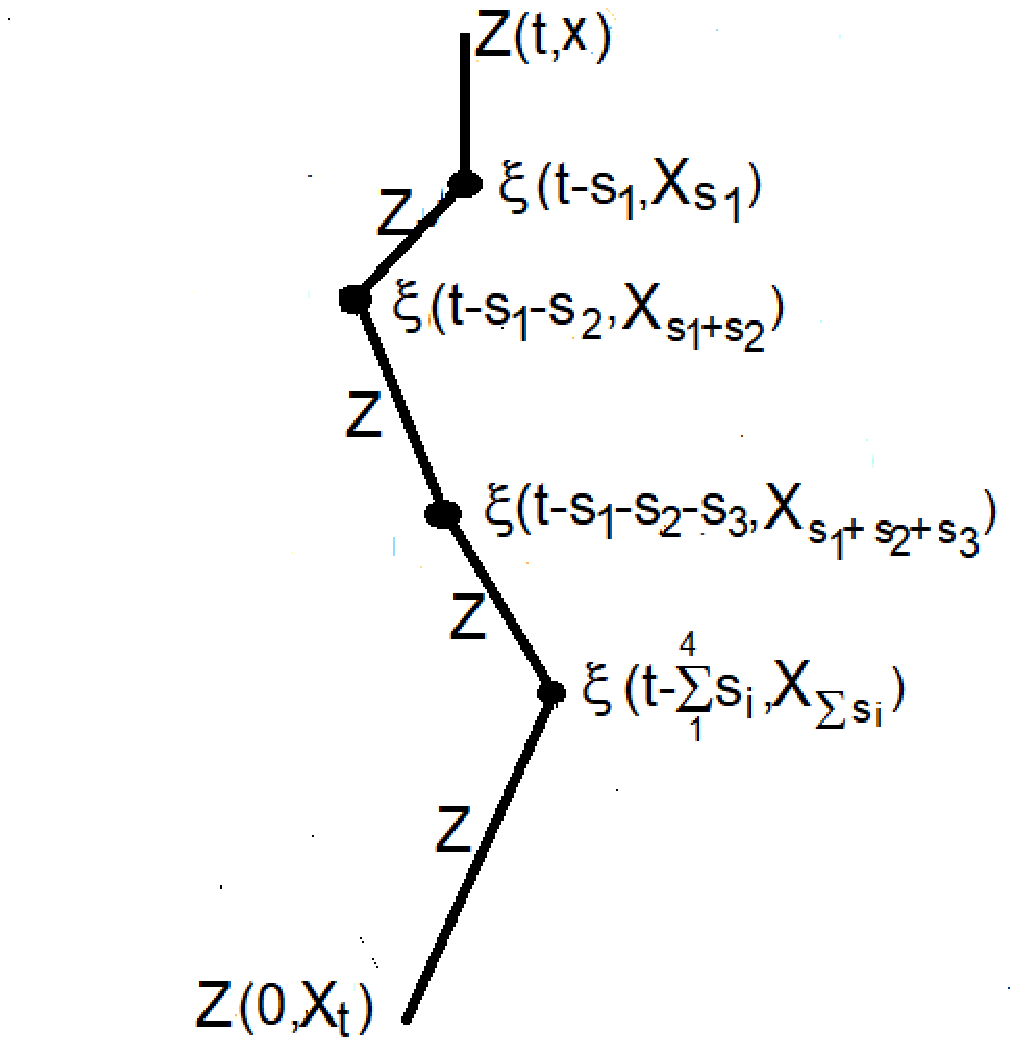}
\caption{A sample path of the Z process}
\label{Z-process}
\end{figure}

Hence, noticing that $e^{t\partial _{x}^{2}}$ is the evolution operator for
a diffusion process $X_{t}$ and that $e^{-t}+\int_{0}^{t}e^{-s}ds=1$, the
second equation in (\ref{K3}) may be written as an expectation value over a
branching and diffusion process which starts at $\left( t,x\right) $ and
evolves backwards in time to $t=0$.%
\begin{equation}
Z\left( t,x\right) =\mathbb{E}_{(t,x)}\left\{ e^{-t}Z\left( 0,X_{t}\right)
+\int_{0}^{t}dse^{-s}Z\left( t-s,X_{s}\right) \lambda \xi \left(
t-s,X_{s}\right) \right\}   \label{K4}
\end{equation}%
$p_{t}=e^{-t}$ is the surviving probability from time $t$ to time zero and $%
\rho _{s}=e^{-s}$ is the probability density for branching between $s$ and $%
s+ds$. Between\ branchings the process propagates as a pure diffusion
process and at each branching there is a sampling of the white noise at that
space-time point and the creation of a new propagation path of the process.
A typical sample path is shown in Fig.\ref{Z-process}.

The contribution to the expectation in (\ref{K4}) of a sample path with $n$
branchings is%
\begin{equation*}
F_{n}=\lambda ^{n}\xi \left( t-s_{1},X_{s}\right) \xi \left(
t-s_{1}-s_{2},X_{s_{1}+s_{2}}\right) \cdots \xi \left( t-s_{1}-s_{2}\cdots
-s_{n},X_{s_{1}+s_{2}\cdots +s_{n}}\right) \cdots Z\left( 0,X_{t}\right)
\end{equation*}

The stochastic construction of the solution to equation (\ref{K2}) involves
two distinct probability spaces, $\Omega $ and $\Omega ^{\prime }$, the
first being the auxiliary probability space of the Brownian motion $X_{t}$
used to compute the expectation $\mathbb{E}_{\left( t,x\right) }$ and the
second the probability space of the white noise that drives the equation.
The expectation value in Eq.(\ref{K4}) is over the branching and diffusion
process, considering a particular fixed realization of the white noise $%
\left\{ \xi \left( t,x\right) ,t\in \left[ 0,t\right] \right\} $. Because
the driving term is a white noise, $Z\left( t,x\right) $ in Eq.(\ref{K4}) is
also a random variable in $\Omega ^{\prime }$. If one lets the $\xi $ white
noise realization change at each branching one obtains a process $Z^{\prime
}(t,x)$ in $\Omega \otimes \Omega ^{\prime }$. The solution to (\ref{K2}) is
a partial average over $\Omega $.

So far one has avoided the product of distributions arising from the last
term in Eq.(\ref{K2}) by trading off products by branchings. Then, to
compute the $F_{n}$ functional of each sample path, one deals with the
product of $n$ white noises which, with probability one, are defined at
different space-time points. Therefore one deals with products of $n$
independent Gaussian variables, which are well-defined random variables with
distribution given by Meijer G-functions \cite{produgauss1} \cite%
{produgauss2}. Then, with probability one, the $Z^{\prime }$ process and the
functionals $F_{n}$ are well-defined.

\begin{proposition}
The functionals $F_{n}$ over the branching, diffusion and white noise sample
paths generate a well-defined process $Z^{\prime }$ in $\Omega \otimes
\Omega ^{\prime }$. A partial expectation over $\Omega $, when it exists, is
the $\Omega ^{\prime }-$process solution of Eq.(\ref{K2}). The process
solution of the KPZ equation is obtained from the second equation in (\ref%
{K1}).
\end{proposition}

Here, the second equation in (\ref{K3}) was used to generate a backwards in
time process, from time $t$ to time zero. Similarly, using the first
equation in (\ref{K3}) one obtains a forward in time process, for which a
similar construction may be performed.

The construction of the $Z^{\prime }$ process associated to Eq.(\ref{K2}) is
an exact result. In this case however, in virtue of the simple structure of
the branching process, an heuristic derivation of the expectation value may
be performed. Assuming a single realization of the $\xi $ white noise, the
expectation value over the branching may be performed explicitly by
multiplying each factor by its branching probability $%
e^{-s}e^{-s_{1}}e^{-s_{2}}\cdots e^{-\left( t-\sum s_{i}\right) }$, the last
factor being the final propagation to time zero, and integrating over all
intermediate times. The result is%
\begin{eqnarray}
Z\left( t,x\right)  &=&\mathbb{E}_{\left( t,x\right) }\left\{ e^{-t}Z\left(
0,X_{t}\right) 1+\lambda \int_{0}^{t}dB\left( t-s,X_{s}\right) +\lambda
^{2}\int_{0}^{t}\int_{0}^{s}dB\left( t-s,X_{s}\right) dB\left( t-s-s^{\prime
},X_{s+s^{\prime }}\right) \right.   \notag \\
&&\left. +\lambda ^{3}\int_{0}^{t}\int_{0}^{s}\int_{0}^{s^{\prime }}dB\left(
t-s,X_{s}\right) dB\left( t-s-s^{\prime },X_{s+s^{\prime }}\right) dB\left(
t-s-s^{\prime }-s^{\prime \prime },X_{s+s^{\prime }+s^{\prime \prime
}}\right) +\cdots \right\}   \notag \\
&=&\mathbb{E}_{\left( t,x\right) }\left\{ e^{-t}Z\left( 0,X_{t}\right)
e^{\lambda \int_{0}^{t}dB\left( t-s,X_{s}\right) }\right\}   \label{K6}
\end{eqnarray}%
$B\left( t,x\right) $ being the Brownian motion of which the white noise $%
\xi \left( t,x\right) $ is the distributional derivative. Finally%
\begin{equation}
h\left( t,x\right) =\frac{1}{\lambda }\log \mathbb{E}_{\left( t,x\right)
}\left\{ e^{\lambda h\left( 0,X_{t}\right) }e^{\lambda \int_{0}^{t}dB\left(
t-s,X_{s}\right) }\right\}   \label{K7}
\end{equation}%
which in the $\lambda \rightarrow 0$ limit converges to the known solution
of the stochastic heat equation%
\begin{equation}
h_{\lambda =0}\left( t,x\right) =\mathbb{E}_{\left( t,x\right) }\left\{
h\left( 0,X_{t}\right) +\int_{0}^{t}dB\left( t-s,X_{s}\right) \right\} 
\label{K8}
\end{equation}%
$e^{\lambda \int_{0}^{t}dB\left( t-s,X_{s}\right) }$ and $%
\int_{0}^{t}dB\left( t-s,X_{s}\right) $ are random variables in $\Omega
^{\prime }$ for each sample path of the process in $\Omega $.

So far, as in the work of other authors \cite{Bertini} \cite{Sasamoto} \cite%
{Balazs}, a version of the Cole-Hopf transformation was used to solve the
KPZ equation. However the stochastic solution technique may also be applied
directly to the equation (\ref{1.1}). Because in that equation there is no
natural branching clock, one adds and subtracts a linear term to the equation%
\begin{equation}
\partial _{t}h\left( t,x\right) =\partial _{x}^{2}h\left( t,x\right) -\mu
h\left( t,x\right) +\lambda \left( \partial _{x}h\right) ^{2}+\mu h\left(
t,x\right) -\xi \left( t,x\right)  \label{K9}
\end{equation}%
which allows to write the integral equation as%
\begin{equation}
h\left( t,x\right) =e^{-t\mu }e^{t\partial _{x}^{2}}h\left( 0,x\right)
+\int_{0}^{t}dse^{-\mu s}e^{s\partial _{x}^{2}}\left\{ \lambda \left(
\partial _{x}h\left( t-s,x\right) \right) ^{2}+\mu h\left( t-s,x\right) -\xi
\left( t-s,x\right) \right\}  \label{K10}
\end{equation}%
The simultaneous addition and subtraction of the linear term may seem an
arbitrary procedure. However, as will be seen in the next section, it is a
most natural procedure whenever one wants to provide a probabilistic
interpretation to this type of equation. In this case it plays the same role
as the multiplication factor $e^{-t}$ in the Cole-Hopf transformation of Eq.(%
\ref{K1}).

\begin{figure}[htb]
\centering
\includegraphics[width=0.8\textwidth]{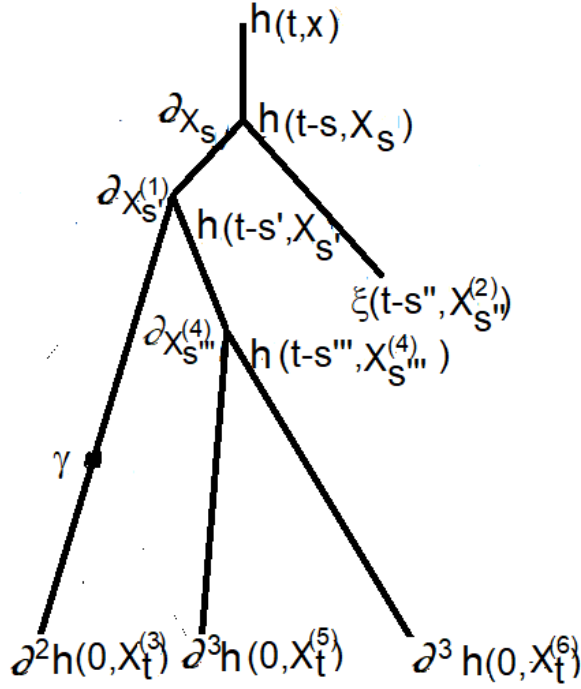}
\caption{A sample path of the DB1-process}
\label{DB1-process}
\end{figure}

The solution will then be an expectation value of a diffusion and branching
process, denoted $DB1$%
\begin{eqnarray}
h\left( t,x\right) &=&\mathbb{E}^{DB1}\left\{ e^{-t\mu }h\left(
0,X_{t}\right) +\gamma \int_{0}^{t}ds\mu e^{-\mu s}\left[ \frac{\lambda }{%
\gamma }\left( \partial _{X_{s}}h\left( t-s,X_{s}\right) \right) ^{2}\right.
\right.  \notag \\
&&\left. \left. +\frac{\mu }{\gamma }h\left( t-s,X_{s}\right) -\frac{1}{%
\gamma }\xi \left( t-s,X_{s}\right) \right] \right\}  \label{K11}
\end{eqnarray}%
with $\gamma =\frac{\lambda +\mu +1}{\mu }$. The $DB1-$process starts from $%
\left( t,x\right) $ and diffuses backwards-in-time either to time zero with
probability $e^{-t\mu }$ or to a branching point at time $s$ with
probability $ds\mu e^{-\mu s}$. At the branching point, with probability $%
\frac{1}{\gamma }$ the process samples the white noise, with probability $%
\frac{\mu }{\gamma }$ the process proceeds undisturbed and with probability $%
\frac{\lambda }{\gamma }$ two new $DB1-$processes are started from $\left(
t-s,X_{s}\right) $ and a derivative label $\partial _{X_{s}}$ is assigned to
the branching point. In the subsequent processes both $h\left(
t-s,X_{s}\right) $ and the label $\partial _{X_{s}}$ are transported by the
processes%
\begin{equation*}
e^{s\partial _{x}^{2}}\partial _{x}h\left( t-s,x\right) =e^{s\partial
_{x}^{2}}\partial _{x}e^{-s\partial _{x}^{2}}e^{s\partial _{x}^{2}}h\left(
t-s,x\right)
\end{equation*}%
Because the propagation process here is a simple diffusion the transport of
the labels is very straightforward and at the final step, when the processes
reach time zero, the derivative is applied at the same point as the point
that is reached by the processes. For more complex propagation processes the
situation is more delicate. This is discussed in the Appendix. To use this
method of labelled branching trees allows for a simple construction of
stochastic solutions for equations involving derivatives or even
nonpolynomial interactions. This method was first used in \cite{VilelaJMP}
(see also \cite{VilelaExt}), being later rediscovered by other authors \cite%
{Touzi}. Fig.\ref{DB1-process} shows one sample path of the $DB1-$process.

The contribution of this sample path to the expectation in Eq.(\ref{K11})
would be%
\begin{equation*}
-\gamma ^{4}\partial ^{2}h\left( 0,X_{t}^{(3)}\right) \partial ^{3}h\left(
0,X_{t}^{(5)}\right) \partial ^{3}h\left( 0,X_{t}^{(6)}\right) \xi \left(
t-s^{\prime \prime },X_{s^{\prime \prime }}^{(2)}\right)
\end{equation*}%
Notice that at each branching point there is a coupling constant $\gamma $.
Also when more than one white noise contribution appears in the
multiplicative functional, with probability one they sample the white noises
at different space-time points. Therefore one has a well-defined product of
independent Gaussian random variables.

As before one deals with two probability spaces $\Omega $ and $\Omega
^{\prime }$. The expectation value in Eq.(\ref{K11}) is an expectation value
over the $DB1-$process for each fixed realization of the white noise
process, hence an expectation in $\Omega $. When a different realization of
the white noise is sampled whenever it appear in the branching tree one
obtains a $DB1^{\prime }-$ process in $\Omega \otimes \Omega ^{\prime }$.

\begin{proposition}
The partial expectation (in $\Omega $) of the multiplicative functional of
the $DB1^{\prime }-$process in $\Omega \otimes \Omega ^{\prime }$(diffusion,
branching and independent white noises at each branching) generates a
solution of the KPZ equation (\ref{1.1})
\end{proposition}

\section{A near mean field phase transition dynamics equation}

Equation (\ref{1.2}) written as an integral equation is

\begin{equation}
\Phi \left( t,x\right) =e^{t\Delta }\Phi \left( 0,x\right)
-\int_{0}^{t}e^{s\Delta }\left( \Phi ^{3}-\xi \right) \left( t-s,x\right) ds
\label{MF1}
\end{equation}%
which may be rewritten as an expectation over the Brownian process $X_{t}$
generated by the Laplacian $\Delta $,%
\begin{equation}
\Phi \left( t,x\right) =\mathbb{E}_{\left( t,x\right) }\left\{ \Phi \left(
0,X_{t}\right) -\int_{0}^{t}\left( \Phi ^{3}-\xi \right) \left(
t-s,X_{s}\right) ds\right\}  \label{MF2}
\end{equation}%
As has been seen in the KPZ example, the replacing of the ill-defined
distributional products by branchings relies on the existence of a
representation of the equation as a branching process. In the Cole-Hopf
transformed KPZ equation (\ref{K2}), the branching clock is supplied by the
linear term $-Z$. In the equation (\ref{MF2}) there is no such term and
multiplying $\Phi $ by $e^{-t}$ is also not convenient because of the $\Phi
^{3}$ term. There is a rather trivial way to obtain a branching clock.
Nevertheless it is useful to explain why, although trivial, it is natural
and rather general. Take a very general approach: One wants to find a
functional of $\Phi $, $\mathcal{F}_{\beta }\left( \Phi \right) $, depending
on a parameter $\beta $, which has a diffusion and branching representation
and is such that in some limit ($\beta \rightarrow 0$, for example) obeys
the equation (\ref{MF2}) and $\lim_{\beta \rightarrow 0}\mathcal{F}_{\beta
}\left( \Phi \right) \rightarrow \Phi $. By hypothesis $\mathcal{F}_{\beta
}\left( \Phi \right) $ obeys a branching and diffusion equation. Therefore%
\begin{eqnarray}
\mathcal{F}_{\beta }\left( \Phi \right) \left( t,x\right) &=&e^{-\lambda
t}e^{t\Delta }\mathcal{F}_{\beta }\left( \Phi \right) \left( 0,x\right)
-\int_{0}^{t}\lambda e^{-\lambda s}e^{s\Delta }\Psi _{\beta }\left( x,\Phi
\left( t-s,x\right) \right) ds  \notag \\
&=&\mathbb{E}_{\left( t,x\right) }\left\{ e^{-\lambda t}\mathcal{F}_{\beta
}\left( \Phi \right) \left( 0,X_{t}\right) -\int_{0}^{t}\lambda e^{-\lambda
s}\Psi _{\beta }\left( X_{t},\Phi \left( t-s,X_{t}\right) \right) ds\right\}
\notag \\
&&  \label{MF3}
\end{eqnarray}%
where $\Psi _{\beta }\left( X_{t},\Phi \left( t-s,X_{t}\right) \right) $ is
the branching kernel.

Compute%
\begin{eqnarray}
&&\mathbb{E}_{\left( t,x\right) }\left\{ \int_{0}^{t}\lambda \mathcal{F}%
_{\beta }\left( \Phi \right) \left( t-s,X_{s}\right) ds\right\}  \notag \\
&=&\mathbb{E}_{\left( t,x\right) }\left\{ \int_{0}^{t}\lambda e^{-\lambda
\left( t-s\right) }\mathcal{F}_{\beta }\left( \Phi \right) \left(
0,X_{s+t-s}\right) ds\right.  \notag \\
&&\left. -\int_{0}^{t}\lambda ds\int_{0}^{t-s}\lambda e^{-\lambda s^{\prime
}}\Psi _{\beta }\left( X_{s+s^{\prime }},\Phi \left( t-s-s^{\prime
},X_{s+s^{\prime }}\right) \right) ds^{\prime }\right\}  \label{MF4}
\end{eqnarray}%
Now summing Eqs.(\ref{MF3}) and (\ref{MF4}) one obtains%
\begin{equation}
\mathcal{F}_{\beta }\left( \Phi \right) \left( t,x\right) =\mathbb{E}%
_{\left( t,x\right) }\left\{ \mathcal{F}_{\beta }\left( \Phi \right) \left(
t,X_{t}\right) -\lambda \int_{0}^{t}ds\left( \Psi _{\beta }\left( X_{s},\Phi
\left( t-s,X_{s}\right) \right) -\mathcal{F}_{\beta }\left( \Phi \right)
\left( t-s,X_{s}\right) \right) \right\}  \label{MF6}
\end{equation}%
This is an useful identity which relies on the Markov nature of the process
(for a detailed proof refer to \cite{Dynkin1} or \cite{Vilela-versus}).

Comparing now (\ref{MF2}) with (\ref{MF6}) and requiring $\mathcal{F}_{\beta
}\left( \Phi \right) \underset{\beta \rightarrow 0}{\rightarrow }\Phi $ one
concludes%
\begin{equation}
\lambda \left( \Psi _{\beta }\left( X_{s},\Phi \left( t-s,X_{s}\right)
\right) -\mathcal{F}_{\beta }\left( \Phi \right) \left( t-s,X_{s}\right)
\right) \underset{\beta \rightarrow 0}{\rightarrow }\left( \Phi ^{3}-\xi
\right) \left( t-s,X_{s}\right)  \label{MF7}
\end{equation}%
that is, the branching kernel $\Psi _{\beta }$ must be of the form%
\begin{equation}
\frac{1}{\lambda }\left( \lambda \Phi +\Phi ^{3}-\xi \right)  \label{MF8}
\end{equation}%
It must contain a branching to $\Phi ^{3}$, a sampling of the white noise $%
\xi $ and a linear term $\lambda \Phi $. The linear term simply means that
when the branching clock rings, with probability $\frac{\lambda }{\lambda +2}
$ the path is not interrupted. The two other possibilities have
probabilities $\frac{1}{\lambda +2}$. The minus sign is a coupling constant
that appears in the contribution to the expectation value when the white
noise is sampled. Alternatively, in the language of superprocesses, it means
that one is dealing not with a superprocess in measures, but with a
superprocess in signed measures \cite{VilelaSPOU}.

\begin{figure}[htb]
\centering
\includegraphics[width=0.8\textwidth]{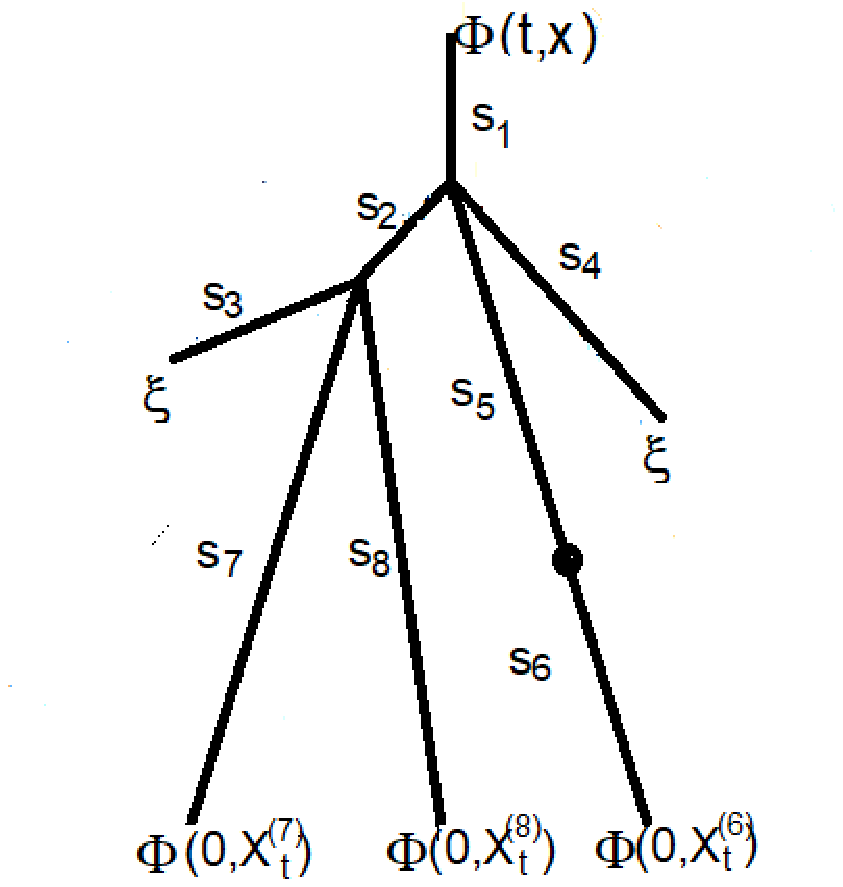}
\caption{A sample path of the DB2-process}
\label{DB2-process}
\end{figure}

Coming back to the original equation one sees that it all amounts to add and
subtract a term linear in $\Phi $ to the equation. One term generates the
branching clock, the other goes to the branching kernel. The result is%
\begin{eqnarray}
\Phi \left( t,x\right)  &=&e^{-t}e^{t\Delta }\Phi \left( 0,x\right)
-\int_{0}^{t}e^{-s}e^{s\Delta }\left( \Phi ^{3}+\Phi -\xi \right) \left(
t-s,x\right) ds  \notag \\
&=&\mathbb{E}_{\left( t,x\right) }^{BD2}\left\{ e^{-t}\Phi \left(
0,X_{t}\right) -3\int_{0}^{t}e^{-s}\frac{1}{3}\left( \Phi ^{3}+\Phi -\xi
\right) \left( t-s,X_{s}\right) ds\right\}   \notag \\
&&  \label{MF9}
\end{eqnarray}%
where again the equations should be interpreted as symbolic definitions of a
process. This backwards-in-time diffusion and branching process ($DB2$),
starting from $\left( t,x\right) $, generates the solution. The probability\
to reach time zero is $e^{-t}$ and the probability density to branch at time 
$s$ is $e^{-s}$. At the branching point there are $3$ possibilities with
equal probability $\frac{1}{3}$. Either the process branches into 3 similar
processes, or the process proceeds undisturbed, or it samples the white
noise. A sample path is shown in Fig.\ref{DB2-process}.

The contribution of this sample path to the expectation in Eq.(\ref{MF9})
would be%
\begin{equation*}
-3^{5}\Phi \left( 0,X_{t}^{(7)}\right) \Phi \left( 0,X_{t}^{(8)}\right) \Phi
\left( 0,X_{t}^{(6)}\right) \xi \left( t-s_{1}-s_{4},X_{s_{1}+s_{4}}\right)
\xi \left( t-s_{1}-s_{2}-s_{3},X_{s_{1}+s_{2}+s_{3}}\right)
\end{equation*}%
As before, when more than one white noise contribution appears in the
multiplicative functional, with probability one they sample the white noises
at different space-time points. Hence a well-defined product of independent
Gaussian random variables is obtained.

The expectation value in Eq.(\ref{MF9}) is an expectation value over the $%
DB2-$process (in $\Omega $) for each fixed realization of the white noise
process. If a different realization of the white noise is sampled whenever
it appear in the branching tree, a $DB2^{\prime }-$ process is obtained in $%
\Omega \otimes \Omega ^{\prime }$.

\begin{proposition}
Partial expectations of the multiplicative functionals associated to the $%
DB2^{\prime }-$process generate the $\Omega ^{\prime }-$process solution of
the equation (\ref{MF1})
\end{proposition}

\section{Conclusions}

The technique of stochastic solutions, by trading products by branchings,
avoids some of the problems associated to the product of distributions.
Because in each sample path the white noises are, with probability one,
sampled at different space-time points, the functionals for each sample path
of the $\Omega \otimes \Omega ^{\prime }-$processes are well-defined.
However the computation and existence of the partial expectation in $\Omega $
may be a delicate matter. In the case of the deterministic equations studied
in the past, existence of the expectation was proved by imposing boundedness
conditions on the multiplicative terms that appear in the functionals. That
procedure is not applicable here because of the unbounded nature of the
white noises. For practical purposes one might impose some truncation on the
stochastic driving term, but this is not satisfactory from a mathematical
point of view.

The computation of the averages may involve an arbitrary large number of
white noise factors. Therefore it is natural to think that the $\Omega
^{\prime }-$process solution will not be distribution -valued, because such
unbounded sums are known not to be in general in the distribution space. A
reasonable conjecture is that the $\Omega ^{\prime }-$process would be
ultradistribution-valued.

\section*{Appendix: The transport of the labels}

For PDE's with derivatives or nonpolynomial interactions in the nonlinear
terms \cite{VilelaJMP} \cite{VilelaExt} the branching and diffusion
processes that generate the solution, receive, at each branching point,
labels that represent the derivatives or the nonpolynomial functions. As the
solution is transported by the diffusion process between the branching
points, so are the labels and, in general, both the solution $h$ and the
labels may be modified by this transport. Let us take as an example a
derivative label that is transported from $\left( x_{0},t_{0}\right) $ to $%
\left( x_{1},t_{1}\right) $ and consider a propagation by a general
diffusion and drift process where both the drift and the diffusion
coefficients may be space-time dependent.%
\begin{equation*}
dW=\overline{b}dt+\overline{\sigma }dB
\end{equation*}%
Then given the value $\partial _{x}h\left( x_{0}\right) $ at $t_{0}$, one
obtains the expectation value of the derivative at $t_{0}+\Delta t$ by
partial integration of the Green's function%
\begin{eqnarray*}
\mathbb{E}\left( \partial _{x}h\left( t_{0}+\Delta t\right) \right) &=&\int
dy\partial _{x}h\left( x_{0}+y\right) \frac{1}{\sqrt{2\pi \overline{\sigma }%
\Delta t}}e^{-\frac{\left( y-\overline{b}\Delta t\right) ^{2}}{2\overline{%
\sigma }\Delta t}} \\
&=&\mathbb{E}\left( h\left( x_{0}+W_{\Delta t}\right) \frac{\left( W_{\Delta
t}-\overline{b}\Delta t\right) }{\overline{\sigma }\Delta t}\right)
\end{eqnarray*}%
Therefore keeping track of the successive increments $W_{\Delta t}$ of the
process one may obtain the value of the derivative at $\left(
x_{1},t_{1}\right) $ from the value of the function that is sampled at $%
\left( x_{0},t_{0}\right) $. This is what some authors \cite{Fournie} have
called \textit{automatic differentiation}. Likewise for the second
derivative one obtains%
\begin{eqnarray*}
\mathbb{E}\left( \partial _{x}^{2}h\left( t_{0}+\Delta t\right) \right)
&=&\int dy\partial _{x}^{2}h\left( x_{0}+y\right) \frac{1}{\sqrt{2\pi 
\overline{\sigma }\Delta t}}e^{-\frac{\left( y-\overline{b}\Delta t\right)
^{2}}{2\overline{\sigma }\Delta t}} \\
&=&\mathbb{E}\left( h\left( x_{0}+W_{\Delta t}\right) \frac{\left( W_{\Delta
t}-\overline{b}\Delta t\right) ^{2}-\overline{\sigma }\Delta t}{\left( 
\overline{\sigma }\Delta t\right) ^{2}}\right)
\end{eqnarray*}

However, whenever $\overline{\sigma }$ and $\overline{b}$ are space
independent constants, one sees from the linear evolution that takes place
between the branching points%
\begin{equation*}
\frac{\partial h}{\partial t}=\frac{1}{2}\overline{\sigma }^{2}\partial
_{x}^{2}h+\overline{b}\partial _{x}h
\end{equation*}%
that all derivatives of $h$ evolve according to the same process as $h$.
Therefore it is sufficient to apply the derivatives to the sampled initial
condition when this is reached in the process, without any need to do the
backtrack of automatic differentiation. However this no longer applies if at
some intermediate step in the branching a dependency on the $x$ coordinate
is introduced.

For the functional labels $f\left( h\right) $,%
\begin{equation*}
\frac{\partial f\left( h\right) }{\partial t}=\frac{1}{2}\overline{\sigma }%
^{2}\partial _{x}^{2}f\left( h\right) +\overline{b}\partial _{x}f\left(
h\right) -\frac{1}{2}\overline{\sigma }^{2}\left( \partial _{x}f\left(
h\right) \right) ^{2}\frac{\partial ^{2}f\left( h\right) }{\partial h^{2}}
\end{equation*}%
hence $f\left( h\right) $ evolves with the same process as $h$ only if $%
f\left( h\right) $ is a linear function of $h$. If it is not, the
calculation of the values at the branching points should be computed with a
backtracking similar to automatic differentiation, namely%
\begin{equation*}
\mathbb{E}\left( f\left( h\left( t_{0}+\Delta t\right) \right) \right) =%
\mathbb{E}\left( f\left( x_{0}+W_{\Delta t}\right) \right)
\end{equation*}%
Therefore, the simplest situation occurs when $\overline{\sigma }$ and $%
\overline{b}$ are chosen to be constants in the linear part of the equation
together with a polynomial approximation for$\ f\left( h\right) $ because
then it can always be decomposed into a product of first order polynomials.

\end{document}